\documentclass[12pt,a4paper]{article}
\usepackage{amssymb}
\usepackage{amsmath,amsfonts,amsthm}

\usepackage{graphicx}
\usepackage{amsmath}

\usepackage{amsfonts}
\usepackage{amsthm,amscd}
\usepackage{graphicx}
\usepackage{amsmath}
\usepackage{amssymb}
\usepackage{amstext}

\newtheorem{theorem}{Theorem}
\newtheorem{corollary}[theorem]{Corollary}
\newtheorem{definition}[theorem]{Definition}

\newtheorem{lemma}[theorem]{Lemma}

\newtheorem{proposition}[theorem]{Proposition}

\title{ On Bertelson-Gromov Dynamical Morse Entropy}

\author{Artur O. Lopes \, and Marcos Sebastiani}

\begin{document}

\maketitle

\begin{abstract}

In this mainly expository paper
we present a detailed proof of several results contained in a paper by M. Bertelson and M. Gromov on Dynamical Morse Entropy. This is an introduction to the ideas presented in that work.

Suppose $M$ is compact oriented connected $C^\infty$ manifold of finite dimension. Assume that $f_0 :M \to [0,1]$ is a surjective Morse function.

For a given natural number $n$, consider the set $M^n$  and for $x=(x_0,x_1,...,x_{n-1}) \in M^n$, denote
$ f_n (x) = \frac{1}{n} \, \sum_{j=0}^{n-1} f_0 (x_j).$

The Dynamical Morse Entropy describes for a fixed interval $I\subset [0,1]$ the asymptotic growth of the number of critical points of $f_n$ in $I$, when $n \to \infty$.

The part related to the Betti number entropy does not requires the differentiable structure.

One can describe generic properties of potentials defined in the $XY$ model of Statistical Mechanics with this machinery.

\end{abstract}

\section{Introduction}

We follow the main guidelines and notation of \cite{BG}.

A Morse function is a smooth function such all critical points are not degenerate (see \cite{Mil}).


Suppose $M$ is compact oriented $C^\infty$ manifold of dimension $q \geq 1.$  Assume that $f_0 :M \to [0,1]$ is a surjective Morse function and $\Gamma$ is a free group with basis $\gamma_1,...,\gamma_n$. We assume that $f_0$ has $p$ critical points ($p\geq 2$).

Suppose $\Omega \subset \Gamma$ is a finite non-empty set. If $x \in M^\Omega$ we denote $x_\gamma \in M$, $\gamma \in \Omega$, the corresponding coordinate.

Then, we define
$f_\Omega : M^\Omega \to [0,1]$ by the expression
$$ f_\Omega (x) = \frac{1}{| \Omega|} \, \sum_{\gamma \in \Omega} f_0 (x_\gamma),$$
where $|\Omega|$ is the cardinality of $\Omega$. This function $f_\Omega$ is also a surjective Morse function.

\section{The $X\,Y$ model}

As a particular case we can consider $\Gamma=\mathbb{Z}$, the set $M^\mathbb{Z}$  and for $x=(x_j)_{j \in \mathbb{Z}} \in M^\mathbb{Z}$, $n>0$, $f_0: M \to \mathbb{R}$, and
$$ f_n (x) = - \frac{1}{n} \, \sum_{j=0}^{n-1} f_0 (x_j).$$

The question about the minus sign in front of the sum is not important but if we want that $f_0$ represents a kind of energy we will keep the $-$ (at least in this section).

In the model it is natural to consider that
adjacent molecules in the lattice interact via a
potential (energy) which is described by the smooth function of two variables $f_0$. The mean energy up to position $n$ is described by $f_n$. The points $x\in M^n$ where the mean $n$-energy is lower or higher are of special importance.
We are interested here, among other things,  in the growth of the number of critical values, when $n \to \infty$.
The critical points are called the stationary states (see \cite{BG}).

Denote by $\text{Cri}_n (I)$ the number of critical points of $f_n$ in a certain interval $f^{-1}(I)$.
Roughly speaking the purpose of \cite{BG} is to provide for a fixed value $c\in [0,1]$ a topological lower bound for
$$ \lim_{\delta \to 0} \lim_{n \to \infty} \frac{\log (\text{Cri}_n (I)  )}{n},\,\,\text{where}\,\, I=(c-\delta,c+\delta) ,$$
in terms of a certain strictly positive  concave function (a special kind of entropy). This is done by taking into account  the homological behavior of the functions $f_n$.

The so called classical $XY$ model consider the case where $M=S^1$ (see for instance \cite{BCLMS}, \cite{LMMS}, \cite{LMST}, \cite{juan}, \cite{Th}, \cite{FHo}, \cite{CG} or \cite{van}). A function $A: (S^1)^\mathbb{Z} \to \mathbb{R}$ describes interaction between sites on the lattice $\mathbb{Z}$ where the spins are on $S^1$. One is interested in equilibrium probabilities $\hat{\mu}$ on $(S^1)^\mathbb{Z}$ which are invariant for the shift $\hat{\sigma}:(S^1)^\mathbb{Z} \to (S^1)^\mathbb{Z}.$
A point $x$ on $(S^1)^\mathbb{Z}$ is denoted by $x=(...,x_{-2},x_{-1}\,|\,x_0,x_1,x_2,...).$

In the case the potential $A$ depend just on the first coordinate $x_0\in S^1$, that is $A(x)=f_0(x_0)$, then the setting described above applies.

In the case the potential $A$ depend just on the two first coordinate $x_0,x_1\in S^1$, that is $A(x)=f_0(x_0,x_1)$, then, we claim  that the setting described above in the introduction applies. This is the case when $f_0: S^1 \times S^1 \to \mathbb{R}$.
Indeed, in this case  one can take $M = S^1 \times S^1$ and consider that $f_0$ acts on $M$. In this case we can say that $f_0$ depends just in the first coordinate on $M^\mathbb{Z}=(S^1 \times S^1)^\mathbb{Z}$ and adapt  the general formalism we describe here.

Therefore, we will state all results for $f_0 : M \to \mathbb{R}$, that is, the case the potential on $M^\mathbb{Z}$ depends just on the first coordinate.

In the case $\hat{\mu}$ is ergodic $f_n$ describes Birkhoff means which are $\hat{\mu}$ almost everywhere constant. We are here interested more in the topological and not in the measure theoretical point of view.

In the measure theoretical (or Statistical Mechanics) point of view, if one is interested  in equilibrium states at positive temperature $T=1/\beta$, then, is natural to consider  expressions like $\int \,e^{ \sum_{j=0}^{n-1} -\beta\, f_0 (x_j)}\, dx_0\,dx_1...d_{x_{n-1}}$ (or, when the set of spins is finite:
$\sum e^{ \sum_{j=0}^{n-1} - \beta\, f_0 (x_j)} $) and its normalization  (see \cite{Sa2}, \cite{CL1} and \cite{CL2}) which defines the partition function.

By the other hand if one is interested in the zero temperature case (see for instance \cite{BLL}), then, expressions like $ -\sum_{j=0}^{n-1} f_0 (x_j)$ are the main focus. For instance, if $f_0$ has a unique point of minimum $x^{-} \in S^1$, then $\delta_{(x^{-})^\infty}$ defines the ground state (maximizing probability). In the generic case the function $f_0$ has indeed a unique point of minimum.

Given $f_0: M \times M \to \mathbb{R}$ and $n$ one can also consider periodic conditions. In this case we are interested in sums like
$$ \tilde{f}_n (x) = -\frac{1}{n} \, ( f_0 (x_0)+ f_0 (x_1)+...+ f_0 (x_{n-2})+ f_0 (x_0)),$$
or
$$ -( f_0 (x_0)+ f_0 (x_1)+...+ f_0 (x_{n-2})+ f_0 (x_0)).$$

In the case we want to get Gibbs states via the Thermodynamic Limit (see for instance \cite{CL1} or \cite{Sa2}), given a natural number $n$, we have to look for the  probability $\mu$ on $M^n$ (absolutely continuous with respect to Lebesgue probability)
which maximizes
$$\int \,e^{-\, \sum_{j=0}^{n-1} \beta\, f_0 (x_j)}\, d\, \mu (dx_0,\,dx_1,...,d_{x_{n-1}}),$$
or, at zero temperature the periodic probability $\mu$ on $M^n$
which maximizes
$$-\,\int \sum_{j=0}^{n-1} \, f_0 (x_j)\,\, d\, \mu (dx_0,\,dx_1,...,d_{x_{n-1}}).$$

In the last case when there is a unique point $x^{-}$ of minimum for $f_0$ then for each $n$ the solution $\mu$ is a delta Dirac on $(x^{-})^n$.

One can easily adapt the reasoning of \cite{AFMT} to show that for a generic $f_0$ we get that $\tilde{f}_n$ is a Morse function for all $n$.

When $f_0$ is not generic several pathologies can occur (see for instance \cite{van}, \cite{BCLMS} and \cite{juan}).

Suppose the case when there is a unique point $x^{-}$ of minimum for $f_0$. For each $\beta>0$ and  $n$ denote by $\mu_{n,\beta}$ the absolutely continuous with respect to Lebesgue probability
which maximizes
$$\int \,e^{ - \,\sum_{j=0}^{n-1} \beta\, f_0 (x_j)}\, d\, \mu (dx_0,\,dx_1,...,d_{x_{n-1}}).$$

By the Laplace method (adapting  Proposition 3 in \cite{LMST} or Lemma 4 in \cite{LMMS}) we get that when $\beta \to \infty$ and $n \to \infty$ the probability $\mu_{n,\beta}$ converges to the Dirac delta on $(x^{-})^\infty$. Therefore,
in the generic case this last probability is the ground state (zero temperature limit).

\bigskip

\section{ The general model - the  dynamical Morse entropy}

From now we forget the $-$ sign in front of $f_0$. For instance,  $ f_n (x) = \frac{1}{n} \, \sum_{j=0}^{n-1} f_0 (x_j, x_{j+1}).$
\medskip

Given $c \in [0,1]$ and $\delta>0$, take $N_\Omega (c, \delta)$ the number of critical points of $f_\Omega$  in
$f^{-1}_\Omega [c-\delta, c+ \delta]$.  Note that if $f_0$ has $p$ critical points then $f_\Omega$ has $p^{| \Omega|}$ critical points.

Consider the cylinder sets
$$ \Omega_i \,=\, \{a_1\, \gamma_1+...+ a_n \gamma_n\,;\, |a_j| \leq i,\,1\leq \,j\leq n\,\},\,\, i=1,2,....$$

Denote $N_i (c, \delta)= N_{\Omega_i}(c, \delta)$. Then, of course, $N_i(c, \delta)$ for $c$ fixed decrease with $\delta$.

For a fixed $0\leq c \leq 1$, we denote the entropy by
$$ \epsilon(c) = \lim_{\delta \to 0}  \,(\,\liminf_{i \to + \infty}\, \frac{\log ( N_i(\,c, \delta)\,)}{| \Omega_i|} \,).$$

The above limit exists  is bounded by $\log p$ but in principle could take the value $-\infty$. We call $\epsilon(c)$ the {\bf dynamical Morse entropy} on the value $c$.

In the case $\Gamma=\mathbb{Z}$ as we mentioned before $ \epsilon(c)$ is described by
$$  \epsilon(c) = \lim_{\delta \to 0}  \,(\,\liminf_{n \to + \infty}\, \frac{\log ( \text{number of critical points of} \,\,f_n \,\,\text{in}\,\,f^{-1}_n [c-\delta,c+\delta] \,)}{n} \,) .$$

Later we introduce a function $b(c)$ (see Definition \ref{a14} and also Definition \ref{def5}), which will be a topological invariant of $f_0$. The function $b(c)$ is defined in terms of rank of linear operators and  Cohomology groups.

We will show later that

1) $ 0 \leq b(c) \leq \epsilon(c), $ $\,0\leq c \leq 1$;

2) $b(c)$ is continuous and concave;

3) $b(c)$ is not constant equal to $0$.

Finally,  in the case $M=S^1$ (the unitary circle) and $f_0$ has just two critical points, we show in section 7 that
$$\epsilon(c) = b(c) = - c \log c - (1-c)  \log (1-c). $$

$b(c)$ is sometimes called the {\bf Betti entropy} of $f_0$.

Our definition of $b(c)$ is different from the one in \cite{BG} but we will show later (see section \ref{ult}) that is indeed the same.

A key result in the understanding of the main reasoning of the paper  is Lemma \ref{a13} which claims that for any  Morse function $f$, given $a,b\in \mathbb{R}$, $a<b$, the number of critical points of $f$ in $f^{-1} [a,b]$ is bigger or equal to the dimension of the vector space
$$\,\, \frac{H^* ( f^{-1} (\infty,b)\,)}{H^* (f^{-1} (-\infty, a)\,)} ,$$
where $H^*$ denotes the corresponding cohomology groups  which will be defined in the following paragraphs (see also \cite{Mas} for basic definitions and properties).

$
H^* (X,\mathbb{R})$ denotes the usual cohomology. Note that $H^*$ will have another meaning (see definition \ref{HH}).

\bigskip

\section{Cohomology} \label{co}

Suppose $X$ is a metrizable, compact, oriented topological manifold $C^\infty$ manifold. We will consider the singular homology. Suppose $U \subset X$ is an open set and $a \in H^* (X, \mathbb{R}).$ The meaning of the statement supp $a\, \subset U$ is: there exist an open set $V\subset X$, such that, $X= U \cup V$, and $a|_V =0.$

\medskip

\begin{definition} \label{HH} $H^*_X (U)= \{\,a\in H^* (X,\mathbb{R}) \,:\, \text{supp } \,\,a \subset U\},$ where $U$ is an open subset of $X$. When $X$ is fixed
we denote $H^*_X(U) = H^* (U)$.

\end{definition}

Remember (see for instance \cite{Mas}) that when $U\subset X$ is open we get the exact cohomology sequence:
\begin{equation} \label{ex}  ... \to H^{k-1} (X - U,\mathbb{R})  \to H^k_c (U,\mathbb{R})\to H^k (X,\mathbb{R}) \to H^k (X-U,\mathbb{R}) \to H^{k+1}_c (U,\mathbb{R})  \to...
\end{equation}

where $H_c^* $ denotes the support compact cohomology.

\begin{lemma} \label{a3}
If $U$ is an open set, then
$$ H^* (U) = \text{Im} (\,\,H^*_c (U,\mathbb{R}) \to H^* (X,\mathbb{R})\,\,)=\,\text{Ker}\,(\,\, H^* (X,\mathbb{R}) \to H^* (X-U,\mathbb{R})\,\,) .$$

\end{lemma}

{\bf Proof:} The second equality follows from the fact that the above sequence is exact.

We will prove that
$$ \text{Im} (\,\,H^*_c (U,\mathbb{R}) \to H^* (X,\mathbb{R})\,\,)\subset  H^* (U) \subset\,\text{Ker}\,(\,\, H^* (X,\mathbb{R}) \to H^* (X-U,\mathbb{R})\,\,) .$$

Let $a\in \text{Im} (\,\,H^*_c (U,\mathbb{R}) \to H^* (X,\mathbb{R})\,\,)$. Then, $a$ is represented by a cocycle $\alpha$ with compact support $K\subset U$. Therefore, $a\,|(X-K)=0$.

Defining $V=X-K$ we have that $U \cup V=X$ and $a\,|V=0$. Then, $ a \in H^* (U)$.

Let be $\alpha\in H^* (U)$. Let $V\subset X$ be an open set such that $U \cup V=X$ and $\alpha\,| V=0$.

Since $X-U\subset V$, we have $\alpha\,|(X-U)=0$.

Then, $\alpha \in $ Ker $(\,H^* (X, \mathbb{R}) \to H^* (X-U, \mathbb{R}) \,).$
\qed

\begin{lemma}
\label{a4} If $U$ is an open set then $H^* (U)$ is a graded ideal of the ring of cohomology of $X$.
\end{lemma}

{\bf Proof:} This follows at once from Lemma \ref{a3}.

\qed

\medskip
Now we consider  a continuous function $f:X \to \mathbb{R}.$

\medskip

\begin{definition} \label{def5}
Given $\delta>0$ and $c \in \mathbb{R}$ we define
$$ b_{c, \delta} ' = Dim \, ( \,\frac{H^* (f^{-1} (- \infty, c +\delta)\,)\,) }{H^* (f^{-1} (-\infty, c- \delta)\,  )} \,).$$

\end{definition}

\begin{proposition} \label{a9} Suppose $X$ and $Y$ are metrizable compact, oriented  topological manifolds, moreover take $f:X \to \mathbb{R}$, $g:Y \to \mathbb{R}$ continuous functions. If we define $f \oplus g: X \times Y \to \mathbb{R}$, by $(f \oplus g)(x,y) = f(x) + g(y)$, then,
if $c,c'\in \mathbb{R}$, $\delta, \delta '>0$, we get
\begin{equation}  \label{op}
b_{c,\,\delta} ' (f) \, \, b_{c ',\,\delta '} '(g) \,\leq b_{c+ c', \,\delta+ \delta '} ' ( f \oplus g).
\end{equation}

\end{proposition}

Before  the proof of this import proposition we need two more lemmas.

As it is known (see \cite{Mas}) the cup product $\vee$   defines an isomorphism
$$ \mu: H^* (X,\mathbb{R}) \otimes H^* (Y,\mathbb{R}) \,\to\, H^* (X \times Y,\mathbb{R}).$$

\begin{lemma} \label{a10}
If $U\subset X$ and $V\subset Y$ are open sets, then
$$ \mu (\, H^*_X (U) \otimes H^* (Y,\mathbb{R}) +  H^* (X,\mathbb{R}) \otimes H^*_Y (V)\,)\, =H^*_{X \times Y}(\,( U \times Y) \cup (X \times V)\,) .$$

\end{lemma}

{\bf Proof:}  By Lemma \ref{a3} we get
$$H^*_{X \times Y}(\,( U \times Y) \cup (X \times V)\,) = \text{Ker}\,  (\,H^* (X \times Y,\mathbb{R}) \to H^* ((X-U) \times (Y - V),\mathbb{R} ).$$

Then,
$$ H^*_{X \times Y} (\,(U\times Y) \cup (X \times V)\,) = $$
$$\mu(\,\,   \text{Ker}\,  (\,H^* (X,\mathbb{R}) \otimes H^* ( Y ,\mathbb{R})\, \to H^* (X-U,\mathbb{R}) \otimes H^* (Y - V,\mathbb{R}) \,)\,\,).$$

From simple Linear Algebra arguments the claim follows from Lemma \ref{a3}.

\qed

\begin{lemma} \label{a11} If $U \subset X$ and $V \subset Y$ are open sets then
$$ \mu ( \,H^*_X (U) \otimes H^*_Y (V)\,) = H^*_{X \times Y} (U \times V).$$

\end{lemma}

{\bf Proof:} The $\vee$ product defines a natural isomorphism
$$ H^* (X, X-U,\mathbb{R}) \otimes H^* (Y, Y-V,\mathbb{R}) \, \to H^*(X \times Y, (X \times (Y-V) \cup (X-U) \times Y,\mathbb{R})= $$
$$ H^*(\,X \times Y, (X \times Y ) - (U \times V,\mathbb{R})\,). $$

By Lemma \ref{a3} and the exact relative cohomology sequence we get:

$$ H^*_X ( U) =\text{Im} \,(\, H^* (X, X-U,\mathbb{R}) \,\to H^*(X,\mathbb{R})\,), $$

$$ H^*_Y (V) =\text{Im} \,(\, H^* (Y, Y-V,\mathbb{R}) \,\to H^*(Y,\mathbb{R})\,), $$

and

$$ H^*_{X\times Y}  (U \times V) =\text{Im} \,(\, H^* (X \times Y, (X \times Y) - (U \times V,\mathbb{R}))\,\to\, H^* (X \times Y,\mathbb{R})\,) . $$

From this the claims follows at once.

\qed

\bigskip

Now we will present the
proof of Proposition \ref{a9}.

Take $h= f \oplus g$ and denote

$$ A^{-} = f^{-1} (-\infty, c-\delta), \,\,B^{-} = g^{-1} (-\infty, c'-\delta'), \,\,C^{-} = h^{-1} (-\infty, (c+c')-(\delta+ \delta')) \,\,,$$
and

$$ A^{+} = f^{-1} (-\infty, c+\delta), \,\,B^{+} = g^{-1} (-\infty, c'+\delta'), \,\,C^{+} = h^{-1} (-\infty, (c+c')+(\delta+ \delta')) \,\,.$$

Note that
$$ A^+ \times B^+ \subset C^+ \subset (A^+ \times Y) \cup (X \times B^+)$$

$$ A^{-} \times B^{-} \subset C^{-} \subset (A^{-} \times Y) \cup (X \times B^{-}).$$


Consider the commutative diagram
$$H^* (X,\mathbb{R}) \otimes H^* (Y) \,\,\,\,\,\,\,\,\,\,\,\,\,\,\,\,\,\,\,\,\,\,\,\,\,\,\,\,\to\,(\text{using} \, \mu\,)\,\,\,\,\,\,\,\,\,\,\,\,\,\,\,\,\,\,\,\,\,\,\,\,\,\,\,\,\, H^* (X \times Y,\mathbb{R})$$
$$\,\,\,\,\,\,\,\,\,\,\,\,\cup\,\,\,\,\,\,\,\,\,\,\,\,\,\,\,\,\,\,\,\,\,\,\,\,\,\,\,\,\,\,\,\,\,\,\,\,\,\,\,\,\,\,\,\,\,\,\,\,\,\,\,\,\,\,\,\,\,\,\,\,\,\,\,\,\,\,\, \,\,\,\,\,\,\,\,\,\,\,\,\,\,\,\,\,\,\,\,\,\,\,\,\,\,\,\,\,\,\, \,\,\,\,\,\,\,\,\,\,\,\,\,\,\,\,\,\,\,\,\,\,\,\,\,\,\,\,\,\,\, \,\,\,\,\,\,\,\,\,\,\cup\,\,\,\,\,\,\,\,\,\,\,\,\,\,\,\,$$
$$H^*_X (A^{+}) \otimes H^*_Y (B^{+}) \to H^*_{X \times Y} (C^{+}) \subset H^*_{X \times Y} (\, (A^{+} \times Y)\cup (X \times B^{+})\,)$$
$$ \cup \,\,\,\,\,\,\,\,\,\,\,\,\,\,\,\, \,\,\,\,\,\,\,\,\,\,\,\,\,\,\,\,\,\,\,\,\,\,\,\,\,\,\,\,\,\,\,\,\,\,\,\,\,\,\,\,\,\, \,\,\,\,\,\,\,\,\,\,\,\,\,\,\,\,\,\,\,\,\,\,\,\,\,\, \,\,\,\,\,\,\,\,\,\, \,\,\,\,\,\,\cup$$
$$ H^*_X (A^{+}) \otimes H^*_Y (B^{-}) + H^*_X (A^{-}) \otimes H^*_Y (B^{+})\to  H^*_{X \times Y} (\, (A^{-} \times Y)\cup (X \times B^{-})\,)$$
$$ \,\,\,\,\,\, \cup\,\,\,\,\,\,\,\,\,\, \,\,\,\,\,\,\,\,\,\,\,\,\,\,\,\,\,\,\,\,\, \,\,\,\,\,\,\,\,\,\,\cup$$
$$ H^*_ X (A^{-}) \otimes H^*_Y (B^{-}) \to H^*_{X\times Y}(C^{-}).$$

From this follows the linear transformation
$$  \tilde{\mu}: \,    \, \frac{H^*_X (A^{+}) \otimes H^*_Y (B^{+})}{ H^*_X (A^{-}) \otimes H^*_Y (B^{-})}\, \to
\frac{ H^*_{X\times Y}(C^{+})}{ H^*_{X\times Y}(C^{-})}.$$

By the other hand
$$ (\,H^*_X (A^{+}) \otimes H^*_Y (B^{+})\,\cap \mu^{-1} ( \,H^*_{X\times Y}(C^{-})\,)\subset$$
$$(\, H^*_X (A^{+}) \otimes H^*_Y (B^{+})\,\cap \mu^{-1} ( \,H^*_{X\times Y}(\, (A^{-} \times Y) \cup (X \times B^{-})\,)\,)=$$
$$ (\,H^*_X (A^{+}) \otimes H^*_Y (B^{+})\,) \cap (\, \,H^*_{X} (A^{-}) \otimes H^* (Y,\mathbb{R}) + H^*(X,\mathbb{R}) \otimes H^*_Y (B^{-})\,)=$$
$$ \,H^*_X (A^{-}) \otimes H^*_Y (B^{+})\, +  \,H^*_X (A^{+}) \otimes H^*_Y (B^{-})\, .$$

The first equality above follows from Lemma \ref{a10}; the second follows from Linear Algebra; namely,
if $E_2\subset E_1 \subset E$ and $F_2\subset F_1 \subset F$, then
$$ (E_1 \otimes F_1) \cap (E_2 \otimes F + E \otimes F_2) = E_2 \otimes F_1 + E_1 \otimes F_2.$$

From the above it follows that
$$ \text{Ker}\, \tilde{\mu}\, \subset  \frac{ H^*_X (A^{-}) \otimes H^*_Y (B^{+})\, +  \,H^*_X (A^{+}) \otimes H^*_Y (B^{-}) }{ H^*_X (A^{-}) \otimes H^*_Y (B^{-})} .$$

Therefore,

$$b_{c + c', \delta+ \delta'} ' \,=\, \text{dim}\,  \frac{H^*_{X \times Y} (C^+)  }{H^*_{X \times Y} (C^{-}) }\geq \text{dim}\, (\text{Im}\,\tilde{\mu})\geq $$
$$  \text{dim}\, \frac{ H^*_X (A^{+}) \otimes H^*_Y (B^{+})}{ H^*_X (A^{-}) \otimes H^*_Y (B^{+})\, +  \,H^*_X (A^{+}) \otimes H^*_Y (B^{-}) }  =  $$
$$  \text{dim}\,(\,\frac{  H^*_X (A^{+})}{  H^*_Y (A^{-})}\, \otimes  \frac{\,H^*_Y (B^{+})}{ H^*_Y (B^{-})}\,)\,= b_{c,\delta} '(f) \, \,b_{c', \delta'} '(g). $$

\qed

\section{Critical points} \label{cri}

In what follows  $X$ is a compact, oriented  $C^\infty$ manifold and $f:X \to \mathbb{R}$ is a Morse function.

\begin{lemma}
\label{a1} Suppose $X$ is a compact, oriented $C^\infty$ manifold  and $U\subset X$ is an open set. If $a \in H^* (X, \mathbb{R}),$ then, supp $a\subset U$, if and only if, there exists a closed $C^\infty$ differentiable form $w$ such that supp $w \subset U$, and $a$ is the de Rham cohomological class of $w$.

\end{lemma}
\medskip

{\bf Proof:} If there exists $w \in a$, such that supp $w \subset U$, then
$$ a|_{(X- \text{supp}\,\, w)}=0 \, \text{and}\,\,  U \, \cup \, (X- \text{supp}\,\, w)\,=\,X. $$

If there exists an open set $V\subset X$
such that $U \cup V=X $ and $a|_V =0$, then, there exist a $C^\infty$ form $\eta$ on $V$ such that  $d \eta = w|_V$ where $w \in a$.

Let $W$ be an open set such that $\overline{W}\subset V$ and $W\cup U=X$. Take a  $C^\infty$ function $\varphi: X \to [0,1]$ such that
$\varphi|_{\overline{W}} =1$ and  $\varphi|_{X-K} =0$, where $K$ is compact set such that $\overline{W} \subset K \subset V$. Then, $\varphi\, \eta$
has an extension to $X$ and $(w- d \,(\varphi\, \eta) \,)\in a$. But,
$$ \text{supp}\,\,\,(w- d \,(\varphi\, \eta) \,)\subset X-W \subset U.$$

\qed

\begin{lemma} \label{a13}

Given $a,b\in \mathbb{R}$, $a<b$, then, the number of critical points of $f$ in $f^{-1} [a,b]$ is bigger or equal that
$$ dim \,\, \frac{H^* ( f^{-1} (\infty,b)\,)}{H^* (f^{-1} (-\infty, a)\,)} .$$
\end{lemma}

{\bf Proof:}

Without lost of generality we can assume that $a$ and $b$ are regular values of $f$ (decrease $a$ and increase $b$ a little bit).

Given $c_1<c_2<...<c_m$, the critical values of $f$ in $(a,b)$, take
$$ a=d_0<c_1<d_1<c_2<d_2<...<d_{m-1}< c_m< d_m=b.$$

By proposition \ref{a8} and Lemma \ref{a6}, the number of critical points in $f^{-1} (c_i)$, $i=1=,2,...,m$, is bigger or equal to
$$ dim \,\, \frac{H^* ( f^{-1} (\infty,d_i)\,)}{H^* (f^{-1} (-\infty, d_{i-1})\,)} .$$

Finally consider the filtration

$$  \,\, H^* ( f^{-1} (\infty,a)\,=\,H^* (f^{-1} (-\infty, d_0)\,) \subset \,H^* (f^{-1} (-\infty, d_1)\,) \subset ...$$

$$  \,\subset\, H^* ( f^{-1} (-\infty,d_{m-1})\,\subset\,H^* (f^{-1} (-\infty, d_m)\,) = \,H^* (f^{-1} (-\infty, b)\,).$$

\qed

Now we denote $b_\Omega '(c,\delta)=b_{c,\delta}'(f_\Omega)$ and $b_i '(c,\delta)=b_{\Omega_i}'(c,\delta)$, $0\leq c \leq 1,$ $\delta>0$.

\begin{corollary} \label{c13}
$b_i '(c,\delta)\leq N_i(c, \delta)$ for all $i=1,2,3,...$ and $0\leq c\leq 1$, $\delta>0$.

\end{corollary}

\medskip

Now we define the function $b$ using Proposition \ref{a8} a)

\begin{definition} \label{a14}

$$b(c)= \lim_{\delta \to 0} \, \liminf_{i \to \infty} \,\,\frac{\log ( b_i ' (c, \delta))}{| \,\Omega_i\,|}, \,\,  0\leq c \leq 1 .$$

\end{definition}

We will show that in above definition we can change the $\liminf$ by $\lim$.

\begin{lemma} \label{a15}

$$b(c) \leq \epsilon (c)\leq \log(\,\text{ the number of critical points of }\,\,f_0\,).$$

\end{lemma}

{\bf Proof:}
The first inequality follows from corollary \ref{c13}.  From the definition is easy to see that $\epsilon (c)$ is smaller than $\log$ of the number of critical points of $f_0.$

\qed

We denote $B(\Gamma)$ a family of finite subsets of $\Gamma$ and $B_N (\Gamma)$, $N \in \mathbb{N}$, the family of sets $\Omega \in B(\Gamma)$ such that $|\Omega| >N.$

\begin{proposition} \label{a16}

Suppose $\Omega', \Omega''\in B(\Gamma)$ are disjoint not empty sets. Then,
$$ b_{\Omega \cup \Omega ''} ' (\alpha c_1 + (1- \alpha) c_2,\delta) \geq  b_{\Omega '} ' (c_1,\delta)\, b_{\Omega ''} ' (c_2, \delta),  $$
where $0 \leq c_1,c_2\leq 1$, $\delta>0$ and $\alpha = \frac{| \,\Omega '\,|}{|\,\Omega ' \,| + | \,\Omega ''\,|}.$

\end{proposition}

{\bf Proof:}

By definition
$$ f_{\Omega ' \cup \Omega ''} = \alpha f_{\Omega '} \oplus (1- \alpha) f_{\Omega ''}.$$

By Proposition \ref{a9}, as $\delta = \alpha \, \delta + (1- \alpha) \delta,$ then

$$ b_{\alpha \,c_1 + (1- \alpha)\, c_2,\delta} ' (f _{\Omega '\cup \Omega ''}) \geq  b_{\alpha\, c_1,\alpha\, \delta} ' (\alpha \,f _{\Omega '})\,b_{(1- \alpha) \,c_2,(1-\alpha)\,\delta} ' ((1-\alpha)\, f _{ \Omega ''})= $$
$$b_{c_1 , \delta} ' (f_{\Omega '})\, \, \, b_{c_2, \delta} ' (f_{\Omega ''}) .  $$

\qed

\begin{lemma}
\label{a6}
Suppose the interval $[a,b]$ does no contains  critical values of $f$. Then,

$$H^* ( \,f^{-1} (\,-\infty, a\,)\,)= H^* (\, f^{-1} (\,-\infty, b\,)\,).$$
\end{lemma}

{\bf Proof:} This follows from Lemma \ref{a3} and the fact that  $f^{-1} [b,\,\infty)\,$ is a deformation retract of  $f^{-1} [a,\,\infty\,).$

\qed

\begin{definition} \label{a7} Given $c \in \mathbb{R} $ we define
$$ \tilde{b}_c (f) = \lim_{\delta \to 0} b_{c,\delta} ' (f).$$

\end{definition}

\begin{proposition} \label{a8} For a fixed $c$ we have

\medskip

a) $ b_{c,\delta} ' (f)$ decreases with $\delta$ and $b_{c,\delta} ' (f) = \tilde{b}_c (f) $ for all $\delta$ small enough.

\medskip

b) $\tilde{b}_c (f)=0 $ if $c$ is not a critical value of $f$

\medskip

c) $\tilde{b}_c (f) $ is smaller than the number of  critical points of $f$  in $f^{-1}(c)$

\medskip

d) $\sum_c \, \tilde{b}_c (f)=$ Dim $H^* (X).$

\end{proposition}

{\bf Proof:}

a) follows from the above definitions and Lemma \ref{a6}.


\medskip

b) follows from Lemma \ref{a6}

\medskip

For the proof of c) consider the exact diagram

$$\,\,\,\,\, \,\,\,\,\,\,\,\,\,\,\,\,\,\, \,\,\,\,\,\,\,\,\,\,\,\,\, \,\,\,\,\,\,\,\, \,\,\,\,\,\,\,\,\,\,\,\,\, \,\,\,\,\,\,\,\, \,\,\,\,\,\,\,\,\,\,\,\,\, \,\,\,\,H^* (X,\mathbb{R}) $$
$$\,\,\,\,\,\,\,\, \,\,\,\,\,\,\,\,\,\,\,\, \,\,\,\,\,\,\,\,\, \,\,\,\,\,\,\,\, \,\,\,\,\,\,\,\,\, \,\, \,\,\,\,\,\,\,\, \,\,\,\,\,\,\,\,\, \,\,\,\,\,\,\,\,\,\,\,\,\, \,\,\,\,\downarrow\, r_1\,\,\,\,\,\,\,\, \,\,\,\, \,\,\,\,\,\,\,\,\,\,\, \,\,\,r_2 \searrow $$
$$H^* ( f^{-1}[c-\delta,\infty)) ,f^{-1}(c+\delta,\infty),\mathbb{R}\,) \to H^* [ f^{-1}(c-\delta,\infty),\mathbb{R})\to H^* ( f^{-1}[c+\delta,\infty),\mathbb{R}),$$
where $r_1$ and $r_2$ are the restriction homomorphisms.

By lemma \ref{a3}
$$ H^* ( f^{-1}(-\infty,c+\delta))=\text{Ker}\, r_2\,  \,\,\,\text{and}\,\,\,H^* ( f^{-1}(-\infty, c-\delta))=\, \text{Ker}\, r_1.$$

From this follows that
$$ b_{c,\delta} ' (f) = \,\,\text{Dim}\,\, ( r_1 (\text{Ker} \, (r_2)\,  )\,\,\leq\,\,\text{Dim}\,(\, H^* ( f^{-1}(c-\delta,\infty)) ,f^{-1}(c+\delta,\infty)\,),\mathbb{R}\, ) $$
because the above sequence is exact.

In order to finish the proof we apply Morse Theory (see \cite{Mil})
with $\delta$ small enough.

\medskip

For the proof of d) suppose $c_1<c_2<...<c_m$ are the critical values of of $f$. Now, consider
$$ d_0<c_1<d_1<c_2<d_2<...<d_{m-1} <c_m<d_m.$$

Now, from a) and Lemma \ref{a6} we have
$$ \tilde{b}_{c_i} (f) = \text{Dim}\, (  \frac{H^* ( \,f^{-1}(-\infty, d_i)\,) }{H^* (\,f^{-1}(-\infty,d_{i-1})\,)}\,   ),\,\,\,\, i=1,2,...,m.$$

Finally, note that
$$  0=H^* ( \,f^{-1}(-\infty, d_0))\subset H^* ( \,f^{-1}(-\infty, d_1))\subset...\subset
$$
$$ H^* ( \,f^{-1}(-\infty, d_m)) =H^* (X).$$

\qed

\begin{lemma} \label{a17}

Given $\delta>0$, there exists an integer $N$ such that:
$b_\Omega ' (c, \delta)\geq 1$ for all $c \in [0,1]$ and all $\Omega \in B_N (\Gamma)$. Therefore, $b(c)\geq 0$, for all $0\leq c\leq 1$.

\end{lemma}

Before the proof of lemma \ref{a17} we need two more lemmas.

\begin{lemma} \label{a18} Suppose $X$ is  a compact oriented  $C^\infty$ manifold and $f:X \to \mathbb{R}$ is a Morse function. Then,
for all $\delta>0$
$$ b_{a_1, \delta} ' (f) \geq 1\,\, \,\text{and}\,\,\,\, b_{a_2, \delta} ' (f) \geq 1, $$
where $a_1$ and $a_2$ are respectively the maximum and minimum of $f$.

\end{lemma}

{\bf Proof:}
If $\delta$ is small enough, $f^{-1} (-\infty, a_2 + \delta)$ is the disjoint union of a finite number of open discs and
$f^{-1} (-\infty, a_2 -\delta)= \emptyset$.

If $n$ is the dimension of $X$, then, it follows from Lemma \ref{a3}  that

$$ H^n (X,\mathbb{R}) \subset H^* (f^{-1} (-\infty, a_2 + \delta  )\,)\neq 0
$$

and
$$ H^* (f^{-1} (-\infty, a_2 - \delta  )\,)=0.
$$

Then,
$b_{a_2, \delta} ' (f) \geq 1$, if $\delta>0$ is small enough. Therefore, this claim is also true for any $\delta>0$ by Proposition \ref{a8} a).

In a similar way we have that for small $\delta>0$

$$ H^0 (X,\mathbb{R}) \subset H^* (f^{-1} (-\infty, a_1 + \delta  )\,)
$$

and
$$ H^0 (X,\mathbb{R}) \,\,\text{is not contained}\,\,\, H^* (f^{-1} (-\infty, a_1 - \delta  )\,).
$$

From this the final claim is proved.

\qed

\begin{lemma} \label{a19} Consider $\Omega\in B (\Gamma)$ where  $|\Omega|=m\geq 1$, then,
$ b_\Omega ' (k/m, \,\delta) \geq 1,$ for all $\delta> 0$ and $k=0,1,2,...,m$.

\end{lemma}

{\bf Proof:}
If $k=0$, or $m$, the claim follows from Lemma \ref{a18} with $X=M^\Omega$, $f=f_\Omega$.

Given $0,k,m$, $0<k<m$, take $\Omega= \Omega' \cup \Omega ''$, where $\Omega ', \Omega '' $ are disjoints and $k=|\Omega '|$.

By Proposition  \ref{a16} with $c_1=1$ and $c_2=0$ we get
$$ b_\Omega ' (k/m, \,\delta)\geq  b_{\Omega '} ' (1, \,\delta)\, b_{\Omega ''} ' (0, \,\delta)\geq 1. $$

Yet from last lemma.

\qed

Now we will prove Lemma \ref{a17}.

{\bf Proof:}

Take $N> \frac{2}{\delta}$, $\Omega \in B_N (\Gamma) , |\Omega|=m >N$ and $k$ such that $\frac{k}{m} \leq c < \frac{k+1}{m},$

By definition,
$$ b_{c, \delta} ' (f_\Omega) \geq b_{k/m, \, \delta/2} ' (f_\Omega),$$
since $c-\delta<k/m-\delta/2$ and $c+\delta>k/m + \delta/2$.

Therefore, $b_\Omega ' (c, \delta) \geq b_\Omega ' (k/m,\,\delta/2)\geq 1$ by Lemma \ref{a19}.

\qed

\begin{proposition} \label{a20}
$$ 0 \leq b(c) \leq \epsilon (c) \leq \, \log ( \text{number of critical points of}\, f_0\,), \,\, 0\leq c\leq 1.$$

\end{proposition}

{\bf Proof:}
This follows from Lemma \ref{a15} and Lemma \ref{a17}

\qed

\begin{lemma} \label{a21} Given $c\in [0,1]$ and $\delta>0$, consider a non-empty set $\Omega \in B(\Gamma)$ and $\gamma \in \Gamma$.
Then,
$$ b_{\Omega} ' (c,\delta) = b_{\Omega+ \gamma } '(c,\delta).$$

In the case $\Gamma=\mathbb{Z}$ we have that for any $\Omega=\{1,2,...,k\}$
$$ b_{\Omega} ' (c,\delta) = b_{\hat{\sigma}(\Omega)} '(c,\delta),$$
where $\hat{\sigma}$ is the shift acting on $M^\mathbb{Z}$.

\end{lemma}

{\bf Proof:} For fixed $\gamma$ consider the transformation $ x\in M^\Omega \to y\in M^{\Omega + \gamma}$, such that
$y_w =x_{w-\gamma}$, which is a diffeomorphism which commutes $f_{\Omega+ \gamma}$ with $f_\Omega$.

The result it follows from this fact.

\qed

We will show now that indeed one can change $\liminf$ by $\inf$ in Definition \ref{a14}. In order to do that we need the following proposition which describes a kind of subadditivity.

\begin{proposition} \label{a22}

Given an integer number $N>0$ take $h: B_N (\Gamma) \to \mathbb{R}$, $h\geq 0$, which is invariant by $\Gamma$
and such that
$$ h( \Omega ' \cup \Omega '') \geq h(\Omega ') + h (\Omega ''),$$
if $\Omega ', \,\Omega ''\,\in B_N (\Gamma),$ are disjoint. Then, there exists
$$ \lim_{i\to \infty} \frac{h(\Omega_i)}{|\Omega_i|}\geq 0,\,\,\,(\,\text{finite or}\, +\infty\,).$$

\end{proposition}
 \medskip

From this follows:

\begin{corollary}\label{a23}
For $c\in [0,1]$ and $\delta>0$,

a) there exist the limit
$$ \lim_{i \to \infty} \frac{ \log b_i ' (c, \delta)}{|\,\Omega_i\,|}\,=\, b ' (c,\delta).$$
\medskip

b) $0 \leq b ' (c,\delta) \leq$    $\log ( $ number of critical points of\, $f_0$\,),

\medskip

c) $b(c)= \lim_{\delta \to 0} b' (c,\delta) $

\end{corollary}

{\bf Proof:} The claim a) follows from last proposition applied to $h(\Omega) = \log b_\Omega ' (c, \delta)$, by Lemma \ref{a17}, Proposition \ref{a16} taking $c_1=c_2=c$ and also by Lemma \ref{a21}.

Item b) follows from lemma \ref{a16}  and corollary \ref{c13}.

Item c) follows from item a) and the definition of $b(c)$. 

\qed

\medskip

Before the proof of Proposition \ref{a22} we need two lemmas.

\medskip

\begin{lemma} \label{a24} Given an integer positive number $k$, then for each $i > (3\, k\,+1)$ there exists $\Omega_{k,i} \in B(\Gamma)$ such that: a) $\Omega_{k,i} \subset \Omega_i$; \, b) $\Omega_{k,i} $ is a disjoint union of a finite number of translates of $\,\Omega_k$ ; \,c) $\lim_{i \to \infty}
\frac{|\,\Omega_{k,i}|\, }{|\,\Omega_i\,|} =1$; d) $|\Omega_i| - |\Omega_{k,i}|\geq (2k+1)^n$, where $n$ is the number of generators of $\Gamma$.

\end{lemma}

{\bf Proof:} For the purpose of the proof we can assume that $\Gamma= \underbrace{\mathbb{Z} \oplus \mathbb{Z}\oplus ... \oplus \mathbb{Z}}_n$ and take
$\gamma_1,\gamma_2,..,\gamma_n$ the canonical basis.

Take $m \geq 1$ an integer such that
$$ k + \,m\, (2 \,k + \,1) \leq i < k = \, (m+1) \, (2\, k + 1),$$
and
$$ \Omega_{k,i} = \cup \,\,\{\,\,\Omega_k + (\,j_1 (2 k +1),...,j_n (2 k+1)\,)\,|$$
$$\, -m\leq j_1,..,j_n\leq m, \, (j_1,...,j_n)\neq (0,...,0) \,\,\}.$$

It is easy to see that the sets $ \Omega_{k,i} $ satisfy all the above claims.

\qed

\begin{lemma} \label{a25} Given real numbers $x_i\geq 0$ $i=1,2,3,...$, suppose that for each $k$ and each $\epsilon>0$ there exist $N_{k,\epsilon}$ such that
$$ x_i \geq x_k (1-\epsilon) \,\,\text{if}\,\, \, i\geq N_{k,\epsilon}.$$

Then, there exists $\lim_{i \to \infty} x_i$ (which is finite or $+\infty$).
\end{lemma}

{\bf Proof:} Take $L= \limsup_{i \to \infty} x_i$ and $a\in \mathbb{R}$, $a<L$. Then, there exists $x_k>a$. Therefore, $x_i \geq a$, if $i$ is very large. Then, $ \liminf_{i \to \infty} x_i \geq a$. From this follows the claim.

\qed
\bigskip

Now we will prove Proposition \ref{a22}.

\medskip

{\bf Proof:} Suppose $k$ is such that $(2 k +1)^n>N$. Take $i> 3 \, k +1$, then,
$|\,\Omega_{k,i}\,| \geq (2 k + 1)^n>\,N$ and $|\Omega_i - \Omega_{k,i}| \geq (2 k +1)^n >N$.

Then, $h(\Omega_i) = h ( \,\Omega_{k,i}        \,\cup\, (\,\Omega_i - \Omega_{k,i}  )\,)   \geq h(\Omega_{k,i}).$

Moreover, each translate of $\Omega_k$ has cardinality $(2k+1)^n$. Therefore,
$$ h( \Omega_{k,i} )\,\geq\, \frac{|\,\Omega_{k,i} \,|}{|\, \Omega_k \,|} h(\Omega_k).$$

From this follows that
$$ \frac{ h(\Omega_i)}{|\, \Omega_i\,|}\geq  \frac{ h(\Omega_{k,i})}{|\, \Omega_{k,i}\,|}\,\frac{|\,\Omega_{k,i}\,|}{|\, \Omega_i\,|} \geq  \frac{ h(\Omega_k)}{|\, \Omega_{k}\,|}\,\frac{|\,\Omega_{k,i}\,|}{|\, \Omega_i\,|}, $$
and the claim is a consequence of Lemmas \ref{a24} and \ref{a25}.

\qed

\medskip

The next lemma will be used later

\medskip

\begin{lemma}\label{a26} Under the hypothesis of Proposition \ref{a22}  consider
$$ \Omega_i ' = ( \,\Omega_i + (2 i +1) \, \gamma_1\,)\,\cup \Omega_i,\,\,\,i=1,2,3,...$$

Then,
$$ \lim_{i \to \infty} \frac{h( \Omega_i ')}{|\, \Omega_i '\,|}= \lim_{i \to \infty} \frac{h( \Omega_i )}{|\, \Omega_i \,|}.$$

\end{lemma}

{\bf Proof:} If $i >N$, then $|\Omega_i|>N$. Therefore,
$$ h(\Omega_i ') \geq h (\,\Omega_i + \,(\,2\, i +1\,)\, \gamma_1\,)\, +\, h(\Omega_i) = 2 h(\Omega_i).$$

From this follows
$$ \frac{h(\Omega_i ')}{|\,\Omega_i '\,|}\,\geq\,\frac{h(\Omega_i)}{|\,\Omega_i \,|} .$$

Therefore,

$$ \liminf_{i \to \infty}   \frac{h(\Omega_i ')}{|\,\Omega_i '\,|}\,\geq\,\liminf_{i \to \infty}  \frac{h(\Omega_i)}{|\,\Omega_i \,|} .$$

We assume that
$\Gamma= \underbrace{\mathbb{Z} \oplus \mathbb{Z}\oplus ... \oplus \mathbb{Z}}_n$ and
$\gamma_1,\gamma_2,..,\gamma_n$ is the canonical basis.

Take $k$ such that $(2 k +1)^n >N$. For $i> 5 \,k + 2$, take $m>1$ such that $k + m\, (2 k +1) \leq i \leq k + (m+1)\, (2\, k +1).$

Consider
$$ \Omega_{k,i} ' = \cup \,\,\{\,\,\Omega_k ' + (\,j_1 (2 k +1),...,j_n (2 k+1)\,)\,|\, j_1 \,\text{is even}\,, -m\leq j_1\leq m -1,$$
$$\, -m\leq j_2,..,j_n\leq m, \, (j_1,j_2,...,j_n)\neq (0,...,0) \,\,\}.$$

Then,
$ \Omega_{k,i}' \subset \Omega_i$, and $ \Omega_{k,i} '$ is a finite union of disjoints translates of $\Omega_k '$. Moreover $\lim_{i \to \infty} \frac{|\Omega_{k,i} '|}{|\Omega_i|}=1$,
$$ | \,\Omega_{k,i} '\,|\geq 2\, (2 \, k +1)^n >N\,\,\text{and}\,\,\, |\, \Omega_i - \Omega_{k,i} '\,|\geq 2\, (2 \, k +1)^n >N\,\,.$$

From this follows that

$$h( \Omega_i ) = h( \Omega_{k,i} ' \, \cup\, (\,\Omega_i - \Omega_{k,i} ' )\,)\geq h (\,\Omega_{k,i} '\,), $$

By the other hand, all translate of $\Omega_k '$ has cardinality  bigger than $N$.

Therefore,
$$ h( \Omega_{k,i} ' ) \geq \frac{|\, \Omega_{k,i} '\,|}{|\, \Omega_k '\,|}\, h(\Omega_k ') .$$

Then,

$$ \frac{\, h(\Omega_{i} \,)}{|\, \Omega_i \,|}\,\geq \frac{\, h(\Omega_{k,i}' \,)}{|\, \Omega_i \,|}\,\geq\frac{1}{ |\Omega_i|}\,\frac{\, |\,\Omega_{k,i} '\,|\,h(\Omega_{k} ' \,)}{|\, \Omega_k ' \,|}\,=\frac{|\, \Omega_{k,i} ' \,|}{|\,  \Omega_i\,|}\,\,\frac{\, h(\Omega_{k} ' \,)}{|\, \Omega_k '\,|}\,.$$

Now, for a fixed $k$, taking $i \to \infty$ in the above inequality we get
$$ \lim_{i \to \infty} \frac{\, h(\Omega_{i} \,)}{|\, \Omega_i \,|}\,\geq \frac{\, h(\Omega_{k} ' \,)}{|\, \Omega_k '\,|}\,.$$

From this follows that
$$ \lim_{i \to \infty} \frac{\, h(\Omega_{i} \,)}{|\, \Omega_i \,|}\,\geq \limsup_{k \to \infty}\,\frac{\, h(\Omega_{k} ' \,)}{|\, \Omega_k '\,|}\,.$$

\qed

\section{Properties of $b(c)$} \label{bc}

\begin{lemma} \label{a27}

There exists $c\in[0,1]$ such that
$$ b(c) \geq \log ( \,\text{dim}\, H^* (M,\mathbb{R})\,)>0$$
\end{lemma}

{\bf Proof:}

Note that dim \,$ (\,H^* (M)\,)\geq 2$ because dim $M\geq 1$. Let $q$ be the number of connected components of $M$.

If $|\Omega_i|=m_i$, take $0=t_0<t_1<...<t_{m_i} =1$, a partition of $[0,1]$ in $m_i$ intervals of the same size.
By Lemma \ref{a3}
$$ H^* ( f_{\Omega_i}^{-1} (-\infty, t_{m_i})) = \oplus_{r>0} H^r (M^{\Omega_i},\mathbb{R}),$$

Denote $A_{ij}$ a supplement of $H^* ( f_{\Omega_i}^{-1} (-\infty, t_{j-1}))$ in $H^* ( f_{\Omega_i}^{-1} (-\infty, t_{j}))$, $1\leq j\leq m_{i}$.
Then,

$$ \sum_{j=1}^{m_i} \text{dim} \, A_{ij}  =\, \text{dim} \,  H^* ( f_{\Omega_i}^{-1} (-\infty, t_{m_i})) =\, \text{dim}\, H^* (M^{\Omega_i},\mathbb{R})-q. $$

Therefore, there exists  a certain $A_{i \, j_i}=A_i$, such that,
$$ \text{dim}\, A_i\, \geq\, \frac{(\text{dim} \,H^*  (M,\mathbb{R}))^{m_i} -q}{m_i}.$$

Denote $s_i$ the middle point of $(t_{j_i -1}, t_{j_i}]$ and $\delta_i = \frac{1}{2\, m_i}.$

Then, by definition of $ b_i ' (s_i,\delta_i) = $ dim $ A_i$.

There exists a subsequence $s_{i_k}\to c\in [0,1]$,  when $k \to \infty$.

Given $\delta>0$, there exists a $K>0$ such that
$ \delta_{i_k} < \delta/2$ and $|s_{i_k} - c| < \delta/2$, if $k>K$.

This means $c-\delta< s_{i_k} - \delta_{i_k} $ and $s_{i_k} + \delta_{i_k} < c + \delta.$

From this follows that
$ b_{i_k} ' (c,\delta)\geq  b_{i_k} ' (s_{i_k},\delta_{i_k})=$ dim $ A_{i_k}$.

Finally, we get
$$ \frac{ \log (\,b_{i_k} ' (c,\delta)\,)}{|\, \Omega_{i_k}\,|} \geq \frac{1}{m_{i_k} }\, \log \frac{(\text{dim}\, H^* (M,\mathbb{R}))^{m_{i_k}}-q }{m_{i_k}}. $$

Now, taking limit in $k \to \infty$ in  the above expression we get
$$ b ' (c, \delta) \geq \log ( \text{dim}\, (H^* (M,\mathbb{R})).$$

\qed

\begin{lemma}
\label{a28}  The function $b (c)$ is upper semicontinuous.

\end{lemma}

{\bf Proof:} Suppose $c_k$, $k \in \mathbb{N}$ is a sequence of points in $[0,1]$ such that, $c_k \to c$.

Given $\epsilon>0$,
take $\delta>0$, such that, $b ' (c,\delta) < b(c) + \epsilon$. There exists a $N>0$ such that $|c-c_k|< \delta/2$, if $k \geq N.$ Then,
$ c- \delta< c_k - \delta/2$ and $c_k + \delta/2< c + \delta$, if $k \geq N.$

Then, $ b_i ' (c, \delta) \geq b_i ' (c_k, \delta/2)$, if $k\geq N$, for all $i=1,2,3,...$

From this follows that
$b' (c, \delta) \geq b ' (c_k, \delta/2)$. Therefore,
$$ b(c) + \epsilon > b ' (c, \delta)  \geq b ' (c_k , \delta/2) \geq b(c_k) ,\, \text{if}\,\, k\geq N.$$

Therefore
$$ \limsup_{k \to \infty} b(c_k) \leq b(c) + \epsilon,$$
for any $\epsilon>0$. From this it follows the claim.

\qed

\begin{lemma} \label{a29}

The function $b(c)$ is concave.
\end{lemma}

{\bf Proof:} Consider $0\leq c_1< c_2\leq 1$ and $0\leq t\leq 1$, we will show that
$$ b(\, t\, c_1 + (1-t)\, c_2) \geq t\, b(c_1) + (1-t)\, b (c_2).$$

First we will show the claim for $t=1/2$. Denote $\tilde{\Omega}_i = \Omega_i + (2\, i +1) \gamma_1$ and $\Omega_i ' = \Omega_i \cup \tilde{\Omega}_i$.

By Proposition \ref{a16} and Lemma \ref{a21} we get:
$$ b_{\Omega_i '} '(1/2\, c_1 + \,1/2\, c_2, \, \delta) \geq b_{\Omega_i} ' (c_1, \delta)\, b_{ \tilde{\Omega}_i} ' (c_2, \delta)= b_{i} ' (c_1, \delta) \, b_i ' (c_2, \delta),$$
for all $\delta>0$.

Now, applying Lemma \ref{a26} to $h(\Omega) = \log b_{\Omega} ' ( 1/2\, c_1 + 1/2\, c_2, \delta)$, we get
$b ' ( 1/2\, c_1 + 1/2\, c_2, \delta)\geq  1/2\, b ' (\, c_1, \delta)\,  +\,1/2\,b ' (  \, c_2, \delta)$.

Now, taking $\delta \to 0$, we get $b ( 1/2\, c_1 + 1/2\, c_2)\geq   1/2 b (\, c_1)\,+\,  1/2\, b ( \, c_2)$.

The inequality we have to prove is true for a dense set of values of $t$ in $[0,1]$. Then,  by Lemma  \ref{a28} is true for all $t\in[0,1]$.

\qed

\begin{corollary} The function $b(c)$ is continuous for $c \in [0,]$.

\end{corollary}

{\bf Proof:} This follows from Lemmas \ref{a28}  and \ref{a29}.

\qed

We collect all results we get above in the next theorem.

\begin{theorem} \label{a31}

a) $ 0 \leq b(c) \leq  \epsilon (c) \leq \log ($ number of critical points of $f_0$), for all $0 \leq c \leq 1$.

b) $b(c)$ is continuous on $[0,1]$

c) $b(c)$ is concave, that is, its graph is always above the cord

d) $b(c)$ is not constant equal zero. Moreover, there exists a point $c$ where $b(c)\geq \log$ (\,dim $H^* (M,\mathbb{R})\,)>0$
\end{theorem}

\bigskip

\section{An example}

The next example shows that the item d) in the above theorem can not be improved.

\medskip

Take $M=S^n$, $n \geq 1$, and a Morse function $f_0: M \to [0,1]$ which is surjective with only two critical points. Suppose $x_{-}$ is the minimum
and
$x_{+}$ the maximum of $f_0$.
We will compute $b(c)$ and $\epsilon(c)$.

Take  $\Omega \in B(\Gamma)$ with $|\Omega|=m\geq 1$. For each $\Omega' \subset \Omega$ consider the canonical projection $p_{\Omega '}: M^\Omega \to M^{\Omega '}.$ Now, take
$$ \mu^{\Omega '} = p_{\Omega '}^* (\,[\,M^{\Omega '}\,]\, ) \in H^{n\, |\,\Omega '\,|}(M^{\Omega},\mathbb{R}),$$
where $[\,\,\,]$ represents fundamental class. Then,
$$ \{\, \mu^{\Omega '} \, :\, \Omega' \subset \Omega \} $$
is a $\mathbb{R}$-homogeneous basis of $ H^* (M^{\Omega},\mathbb{R}). $

For $0\leq d\leq 1$ denote
$$ L_d =\{ x \in M^\Omega\, : \, f_\Omega (x) <d\, \} \subset M^\Omega.$$

For $x \in M^\Omega$ we denote by $x_\gamma$ the corresponding coordinate, where $\gamma \in \Gamma$.

\begin{lemma} \label{a32} If $0\leq d\leq 1$, where $d$ is not rational, then
$$\{ \mu^{\Omega '} \, :\, |\Omega '|> m\, (1-d)\,\}$$
is a basis of $H^* (L_d)$.

\end{lemma}

{\bf Proof:} Take
$ K_d =  M^\Omega - L_d$. By Lemma \ref{a3}
$$ H^* (L_d) = \text{Ker} (\,H^* (M^\Omega,\mathbb{R}) \,\to \, H^* ( K_d,\mathbb{R})\,)\,\,(\text{natural restriction}). $$

The claim follows from

1) $H^k (M^\Omega,\mathbb{R}) \to H^k (K_d,\mathbb{R}) $ is zero if $k> m \,(1-d)\, n$, and

2) $H^k (M^\Omega,\mathbb{R}) \to H^k (K_d,\mathbb{R}) $ is injective if $k< m \,(1-d)\, n$.

\bigskip

Now we prove (1) and (2).

\medskip

(1) \, Suppose $\Omega' \subset \Omega$ is such that $\mu^{\Omega '} \in H^k (M^\Omega)$ where $k> m\, (1-d)\, n$.
Then, $|\Omega '| > m\, (1-d)$. Suppose
$$ F_{\Omega '} = \{ x\in M^\Omega\, :\, x_\gamma = x_{-}\, ,\,\,\text{if}\, \gamma \in \Omega'\} .$$

If $x \in F_{\Omega '}$, then $f_\Omega (x) \leq \frac{1}{m}\, (m- | \Omega '|)< d$. Then, $F_{\Omega '} \cap K_d=\emptyset.$ This means that:
 if $x \in K_d\, \to\,  x_\gamma \neq x_{-}$ for some $\gamma \in \Omega'.$ Then, $K_d \subset p^{-1}_{\Omega '} (M^{\Omega '} -\,\{z\}\,)$ where $z_\gamma=x_{-}$ for all
 $\gamma \in \Omega'$.

 From this follows
 $$ \mu^{\Omega '}\, |\, K_d\,=\,  p^{*}_{\Omega '}\,(\,[\,M^{\Omega '} \,]\,)\,\,|\,\, K_d\,=0,\,\,\text{because}\,[\,M^{\Omega '} \,]\,\,|\, \, (\,[\,M^{\Omega '} \,]\,-\{z\}\,)=0.$$

\medskip

(2)\, Denote $T=\{x \in M^\Omega \,:\, \text{cardinality} (\, \{\gamma\,:\, x_\gamma = x^{+} \,\})\,>\, m\, d\,\}.$ The set $T$ is closed.

If $x \in T$, then $f_\Omega (x) > \frac{1}{m}\, m\, d =d$. Then, $T\subset K_d$.

We have to show that
$$ H^k (M^\Omega ,\mathbb{R})  \to H^k (T,\mathbb{R})\,\,\text{is injective if}\,\, k<m\, (1-d)\, n.$$

As we had seen before $H^k (M^\Omega,\mathbb{R})=0$ if $k$ is not multiple of $n$. Then, we can assume that $k=q\, n$, if $q=0,1,2,...$. The claim follows from
the next lemma, taking $s$ the integer part of $m\, d$, by the exact sequence of homology, given that $U= U_s (\Omega).$

\begin{lemma} \label{a33} Suppose $s=0,1,2,..,m$. Suppose
$$U_s (\Omega) \,=\, \{x \in M^\Omega \,:\, \text{card}\, ( \,\{\gamma\,:\, x_\gamma= x^{+} \}\,) \,\,\leq s \},$$
then, $H^k_c ( U_s ( \Omega),\mathbb{R}) =0$, if $k <(m-s)\,n$.

\end{lemma}

{\bf Proof:} The claim is trivial for $s=0$ or $s=m$ ($U_0(\Omega)$ is homeomorphic to $(\,\mathbb{R}^n\,)^m$).

The proof is by induction in $m$. The claim for $m=1$ is trivial. Suppose is true for $m-1\geq1$. Take $0<s<m$. Fix $w\in \Omega$ and take $\Omega' = \Omega -\{w\}$.

Consider $\varphi: M^{\Omega '} \to M^\Omega $ and $\psi: M^{\Omega '} \times (M-\{\,x^{+}\,\}) \, \to \, M^\Omega$, where for a given $x$ we  define $\varphi(x)$  by $ x_\omega = x^{+}$ if $ x \in M^{\Omega ' },$ and $\psi (x,u)$ is defined by $x_w =u$ if $x \in M^{\Omega '}$ and $u \in M,$ $u \neq x^{+}.$

$\psi$ identifies $U_s (\Omega ') \times (\,M - \{x^{+}\}\,)$ with an open set $A$ contained in $U_s (\Omega)$.

Moreover, $\varphi$ identifies $U_{s-1} (\Omega ')$ with the complement of this open set $A$ in $U_s (\Omega)$.

As $M - \{x^{+}\}$
is homeomorphic to $\mathbb{R}^n$ and by recurrence we get that

$$ H^k_c ( \,U_s (\Omega ')\, \times\, (\,M - \{x^{+}\}\,)\,,\mathbb{R})=0, $$
if $k< (m -1 -s)\, n +n = (m-s)\, n$ and, moreover, $H^k_c (\,U_{s-1} (\Omega ',\mathbb{R})\,)=0$, if
$ k< (\, (m-1) - (s-1)\,) \, n= (m-s)\, n$.

The exact sequence of homology finish the proof.

\qed

\medskip

Now we fix irrationals
$d_1,d_2$, $0< d_1<d_2<1$. Denote
$  a_m = m\, (1-d_1)$, $ b_m = m (1-d_2),$ and, $ c_m =$ dim $ (\, H^* (L_{d_2})/ H^* (L_ {d_1})\,)$.

By Lemma \ref{a32} we get
$$ c_m = \sum \,\{   \left(\begin{array}{cc}
m\\
j
\end{array}
\right)
:\, b_m < j<a_m \}.$$

Assume $m$ is much more bigger than $(d_2-d_1)$.

Take an integer $j_m$, such that $b_m<j_m<a_m$,
$$    \left(\begin{array}{cc}
m\\
j_m
\end{array}
\right)\, =\, \sup\,  \{   \left(\begin{array}{cc}
m\\
j
\end{array}
\right)
:\, b_m < j<a_m \}.
$$

Then,
$$  \left(\begin{array}{cc}
m\\
j_m
\end{array}
\right)\, \leq c_m \leq (a_m - b_m +1)\, \left(\begin{array}{cc}
m\\
j_m
\end{array}
\right)\,  .$$

By Stirling formula:

$$ \frac{1}{m} \log \left(\begin{array}{cc}
m\\
j
\end{array}
\right)\,  \sim  \frac{1}{m} \, \log (\, \frac{m^{m + \,\,\,1/2}}{j^{j + \,\,\,1/2}\, (m-j)^{\,m-j + \,\,\,1/2} }\,) =$$
$$\frac{1}{m} \log (\,m^{-1/2} \,(\frac{j}{m})^{-1/2} \,( 1- \frac{j}{m})^{-1/2}\, (\frac{j}{m})^{-j}\, (1-\frac{j}{m})^{ -m + j} \,). $$

Therefore,
$$ \frac{1}{m} \log \left(\begin{array}{cc}
m\\
j_m
\end{array}
\right)\,  \sim \frac{1}{m} \log (\, (\frac{j_m}{m})^{-j_m}\, (1-\frac{j_m}{m})^{ -m + j_m} \,)=
$$
$$ - \frac{j_m}{m} \, \log (\frac{j_m}{m})\,-\, (1 - \frac{j_m}{m} )\, \log ( 1- \frac{j_m}{m}),$$
when $m \sim \infty$.

As $1- d_2 < \frac{j_m}{m}< 1- d_1$, then (changing $x$ by $(1-x)$) we get
$$ \limsup_{m \to \infty}
\frac{1}{m} \log
\left(\begin{array}{cc}
m\\
j_m
\end{array}
\right)\,\leq\, \sup_{d_1<x<d_2} \,(\,-\, x \log (x) - (1-x)\,\log (1-x)\,)   ,$$
and
$$ \liminf_{m \to \infty}
\frac{1}{m} \log
\left(\begin{array}{cc}
m\\
j_m
\end{array}
\right)\,\geq\, \inf_{d_1<x<d_2} \,(\,-\, x \log (x) - (1-x)\,\log (1-x)\,)   .$$

From this follows
$$ \limsup_{m \to \infty} \frac{\log c_m}{m}
\,\leq\, \sup_{d_1<x<d_2} \,(\,-\, x \log (x) - (1-x)\,\log (1-x)\,)   ,$$
and
$$ \liminf_{m \to \infty} \frac{\log c_m}{m}
\,\geq\, \inf_{d_1<x<d_2} \,(\,-\, x \log (x) - (1-x)\,\log (1-x)\,)   .$$

\medskip

\begin{proposition} \label{a34}
$$\epsilon (c) = b(c) = -\,c\, \log c -\, (1-c)\, \log(1-c), \,\,0 \leq c\leq 1.$$
\end{proposition}

{\bf Proof:}

Given $ 0<c<1$, there exists small $\delta>0$ such that
$$ 0< c - \delta <c < c+ \delta <1\,\,\,\text{and}\,\,\ c-\delta, c+\delta\,\,\text{are not in }\,\,\mathbb{Q}.$$

From the above for $d_1= c -\delta $ and $d_2= c + \delta$ we get
$$ \inf_{d_1<x<d_2} \,(\,-\, x \log (x) - (1-x)\,\log (1-x)\,)\, \leq b ' (c,  \delta) \leq $$
 $$\sup_{d_1<x<d_2} \,(\,-\, x \log (x) - (1-x)\,\log (1-x)\,)   .$$

Now, taking $\delta \to 0$, we get
$$ b(c)\,=\,(\,-\, c \log (c) - (1-c)\,\log (1-c)\,)   .$$

For $c=0$ or $c=1$ the result follows from continuity.

Now we will estimate $\epsilon(c).$

The critical values of $f_\Omega$ are $0,\frac{1}{m},\frac{2}{m},..., 1$.

To the critical values $\frac{j}{m}$ ($j=0,1,2,..,m$) corresponds
$\left(\begin{array}{cc}
m\\
j
\end{array}
\right)\,$ critical points.

Therefore, given $d_1,d_2 \in \mathbb{R}$ $d_1<d_2$, the number $c_m '$ of critical points of $f_\Omega$ in $ f^{-1}_\Omega (d_1,d_2)$ is

$$ c_m ' = \sum\, \{  \left(\begin{array}{cc}
m\\
j
\end{array}
\right)\, :\, d_1 <\, \frac{j}{m}\, < d_2\, \}=$$
$$\sum\, \{  \left(\begin{array}{cc}
m\\
j
\end{array}
\right)\, :\, m\, (1-d_2) < \,j \,< m\,(1-d_1)\, \} .$$

The computation of $\epsilon(c)$ is analogous to the one for $b(c)$.
This also follows from the last Theorem and the fact that $H^* (M)$ = number
critical points of $f_0$ in the present case.

\qed

\section{About the definition of $b(c)$} \label{ult}

We will show that the definition of $b(c)$ presented here coincides with the one in \cite{BG}.

First we need some preliminary results.

Suppose $X$ is a compact connected oriented $C^\infty$ manifold.

\begin{lemma} \label{a35} Given an open set $V$ in $X$ consider $\alpha\in H^* (X,\mathbb{R})$ such that $\alpha|_V \neq 0.$ Then, there exists $\beta \in H^* (V)$ such that $\alpha \wedge \beta \neq 0.$

\end{lemma}

{\bf Proof:} Take $w \in \alpha$. As $\alpha|_V \neq 0$, then there exists a cycle $z$ on $V$ such that $\int_z w \neq 0.$

Suppose $w'$ is a closed form with compact support on $V$ such that its cohomology class  in $H^*_c (V,\mathbb{R})$ is  the Poincare dual of the homology class of $z$ in $H_* (V,\mathbb{R}).$

$w'$ can be extended to a closed form on $X$ (putting $0$ where needed) and by Poincare duality:
$$ 0 \neq \int_z w\,=\,\int_V w \wedge w'\, = \int_X w\, \wedge w'.$$

Therefore, $w \wedge w'$ is not exact on $X$.

Denote $\beta \in H^* (X,\mathbb{R})$ the cohomology class of $w'$. By Lemma \ref{a3} we have that $\beta \in H^*(V)$. As $w \wedge w'$ is not exact we get that $\alpha \wedge \beta\neq 0.$

\qed
\medskip

Notation: if $S \subset X$, then $\mathcal{H}^* (S)=\cap \,\{ H^* (W) \, : W \subset X$ is an open set and $ S \subset W \}$.

\medskip

\begin{lemma} \label{a36} Suppose $U,V \subset X$ are open sets and $X=U \cup V$. Take $K=U-V$ and $\alpha \in H^* (U).$ Then, $\alpha \wedge \beta =0$ for all $\beta \in H^* (V)$, if and only if, $\alpha \in \mathcal{H}^* (K).$

\end{lemma}

{\bf Proof:} Suppose $\alpha \in \mathcal{H}^* (K)$ and take $\beta \in H^* (V) $. By Lemma \ref{a1} there exists $w \in \beta$ such that supp $w \subset V$.

Take $W= X-$ supp $w$ (which contains $K$).  By definition we get that $\alpha \in H^* (W).$ Then, by Lemma \ref{a1}, there exists $w ' \in \alpha$ such that supp $w' \subset W$. Therefore, $w \wedge w'=0,$ and finally it follows that $\alpha \wedge \beta=0.$

Reciprocally, suppose that $\alpha \wedge \beta =0$ for all $\beta \in H^*(V)$. By Lemma \ref{a35} we have that $\alpha\,|V =0$. Take $W \supset K$, then $V \cup W=X$. Therefore, by definition $\alpha \in H^* (W).$

\qed

\begin{lemma} \label{a37}
Take $K \subset X$ a compact submanifold with boundary  such that $K- \delta K$  is an open subset of $X$.

Then,

$$ \mathcal{H} (K) = \, \text{Ker}\, \,(\,H^* (X,\mathbb{R}) \to H^* (X-K,\mathbb{R})\,)\,\,\text{restriction}.$$

\end{lemma}

{\bf Proof:} Take $W$ an open set by adding a necklace to $K$. Then, $X-K$ can be retracted by deformation over $X-W$.

Then, if $\alpha \in H^* (X,\mathbb{R})$, we get that $\alpha|_{X-K} =0$ is equivalent to $\alpha|_{X-W}=0.$

Now, the claim follows from Lemma \ref{a3} and by the definition of $\mathcal{H}(K).$
\qed

\begin{corollary} \label{a38} Under the same hypothesis of last lemma it also follows that
$ \mathcal{H} (K) = H^* ( $ int  $ (K)\,)$.

\end{corollary}

{\bf Proof:} This follows from the fact that $H^* ( X - $ int $(K)\,,\mathbb{R}) \to H^* (X-K,\mathbb{R})$ is an isomorphism.

\qed

\begin{proposition} \label{a39} Suppose $U,V$ are open sets such that $X= U\cup V$ and moreover that $\overline{U},\overline{V}$ are submanifolds with boundary of $X$.

Consider the linear transformation $L$ such that
$$L: H^* (U) \to \text{Hom}\, (\,H^* (V) , H^* (U \cap V)\,),$$
where, $a \to \,(\, b \to a \wedge b\,).$

Then, the rank of $L$ is dim $ ( \,H^*(U)/ H^* (M - \overline{V})\,).$

\end{proposition}

{\bf Proof:} By Lemma \ref{a36} we get that Ker $L=  H^* (X-V)$. Finally, by the last corollary $H^* (X-V)=  H^* (M-\overline{V}).$

\qed

\medskip

Consider now a Morse function $f:X \to \mathbb{R}$ and $c\in \mathbb{R}$, $\delta>0.$

\begin{definition}\label{a40} $b_{c,\delta} (f)$ is the rank of the linear transformation
$$ H^* (\,f^{-1} (-\infty, c + \delta)\,)\, \to \, \text{Hom}\, (\,H^* (f^{-1} (\, c - \delta, \infty\,)  ,H^* (f^{-1} (c-\delta, c + \delta)\,)\, ) ,$$
where $a \to (b\, \to a \wedge b).$

\end{definition}

\medskip

Note that  $b_{c,\delta} (f)$ decreases with $\delta$.

\medskip

\begin{lemma} \label{a41} If $c-\delta$ and $c+ \delta$ are regular values of $f$, then
$$ b_{c,\delta} (f)= b_{c,\delta}' (f).$$

\end{lemma}

{\bf Proof:} Just apply Proposition \ref{a39} to $U= f^{-1} (-\infty, c + \delta)$ and $V= f^{-1} (\, c - \delta, \infty\,).$

\qed

\medskip

Note that $b_\Omega(c,\delta)= b_{c, \delta} ( f_\Omega)$, where $\Omega \in B(\Gamma)$ and $\Omega \neq \emptyset$, and moreover that $b_i (c,\delta) = b_{\Omega_i} (c, \delta).$ The next limit exists (see \cite{BG}).

\begin{definition} \label{a42}

$$ b(c,\delta) = \lim_{i \to \infty} \, \frac{\log (\,b_i (c, \delta)\,)}{|\, \Omega_i\,|}. $$

\end{definition}

\medskip

The set $S\subset [0,1]$ of all critical values of all $f_\Omega$ is countable. By Lemma \ref{a41} we get that $b_i ' (c, \delta)=b_i (c, \delta)$ if $c-\delta \notin S$  and $c+\delta \notin S$. Therefore, $b'(c,\delta) = b(c,\delta)$ if $c-\delta \notin S$  and $c+\delta \notin S$.

Finally,
$$ \lim_{\delta \to 0}  b ' (c,\delta) =  \lim_{\delta \to 0} b(c,\delta)$$
because both limits exist.

Therefore the function $b(c)$ we define coincides with the one presented in \cite{BG}.

\bigskip

\bigskip

\bigskip

\bigskip

{\bf Instituto de Matematica - UFRGS - Brasil}
\medskip

A. O . Lopes was partially supported by CNPq and INCT.


\bigskip

\begin{thebibliography}{99}




\bibitem{AFMT} M. Asaoka, T. Fukaya, K. Mitsui and M. Tsukamoto, Growth of critical points in one-dimensional lattice systems, preprint Arxiv 2012.


\bibitem{BCLMS}
A. T. Baraviera, L. Cioletti, A. O.  Lopes, J. Mohr and R. R. Souza,
On the general one-dimensional XY model: positive and zero temperature, selection and non-selection. \emph{Rev. Math. Phys.} 23,no. 10, 1063-1113, 82Bxx, 2011.



\bibitem{Bov} A. Bovier.
Statistical Mechanics of Disordered Systems.
A Mathematical Perspective.
Cambridge University Press (2006).


\bibitem{BLL} A. Baraviera, R. Leplaideur and A. O. Lopes, Ergodic Optimization, Zero temperature limits and the Max-Plus Algebra,
mini-course in XXIX Col\'oquio Brasileiro de Matem\'atica - IMPA - Rio de Janeiro (2013)


\bibitem{LMST} A. O. Lopes, J. Mohr, R. R. Souza and P. Thieullen,
Negative Entropy, Zero temperature and stationary Markov chains on the
interval,   \emph{ Bull. Soc. Bras. Math.} Vol 40 n 1,  (2009), 1-52.


\bibitem{LMMS}
A. O. Lopes, J. K. Mengue, J. Mohr and R. R. Souza, Entropy and Variational Principle for one-dimensional Lattice Systems with a general a-priori probability: positive and zero temperature" to appear in \emph{Erg. Theo. and Dyn. Syst}.


\bibitem{BG} M. Bertelson and M. Gromov, Dynamical Morse Entropy,  Modern dynamical systems and applications,
27-44, Cambridge Univ. Press, Cambridge (2004)





\bibitem{Israel} R.B. Israel, Convexity in the theory of lattice gases, Princeton University Press, 1979.



\bibitem{CG}
W. Chou and R. Griffiths, {\it Ground states of one-dimensional systems using effective potentials}, Physical Review B, Vol. 34, N 9, 6219-6234, 1986


\bibitem{juan} D. Coronel and J. Rivera-Letelier, Sensitive dependence of Gibbs measures, preprint (2014)

\bibitem{CL1} L. Cioletti and A. Lopes,
Interactions, Specifications, DLR probabilities and the
Ruelle Operator in the One-Dimensional Lattice

\bibitem{CL2} L. Cioletti and A. Lopes, Phase Transitions in One-dimensional
Translation Invariant Systems:
a Ruelle Operator Approach, to appear in Journ. of Stat. Phys.



\bibitem{FHo} Y. Fukui and M. Horiguchi, One-dimensional Chiral $XY$ Model at finite temperature, Interdisciplinary Information Sciences, Vol 1, 133-149, N. 2 (1995)


\bibitem{FM}
T. Fukaya and M. Tsukamoto, Asymptotic distribution of critical values. Geom. Dedicata 143 (2009), 63–67




\bibitem{Mas}
W. Massey, Homology and Cohomology, M. Dekker (1978)


\bibitem{Mil}
J. Milnor, Morse Theory, Princeton University Press, Princeton (1963)



\bibitem{Si} B. Simon, The Statistical Mechanics of Lattice Gases, Princeton Univ Press, 1993

\bibitem{Sa2} O. Sarig, Lecture notes on thermodynamic formalism for topological Markov shifts, \emph{ Penn State}, 2009.



\bibitem{van} A. C. D. van Enter and W. M. Ruszel.
Chaotic Temperature Dependence
at Zero Temperature. \emph{Journal of Statistical Physics}. Vol. 127. No. 3.
(2007), 567-573.



 \bibitem{Th} C. Thompson.
{\it Infinite-Spin Ising Model in one dimension}. Journal of Mathematical Physics.
(9): N.2  241-245, 1968.

\end{thebibliography}
\end{document}